\documentclass[10pt,twoside,final]{amsart}

\usepackage[english]{babel}
\usepackage{graphicx,epstopdf,epsfig}
\usepackage{amsfonts,epsfig,fancyhdr,graphics,amsmath,amssymb}

\title{Bireflections of the commutator subgroup of an orthogonal group over the reals}

\newtheorem{definition}{Definition}[section]
\newtheorem{theorem}{Theorem}[section]
\newtheorem{proposition}[theorem]{Proposition}
\newtheorem{lemma}[theorem]{Lemma}
\newtheorem{corollary}[theorem]{Corollary}
\newtheorem{remark}[theorem]{Remark}

\newcommand {\ch }{\operatorname{char}}
\newcommand {\SOG }{\operatorname{SO}}
\newcommand {\OG }{\mathrm{O}}

\newcommand {\End }{\mathrm{End}}

\newcommand {\Bahn }{\mathrm{B}}
\newcommand {\Fix }{\mathrm{Fix}}

\newcommand {\Neg }{\mathrm{Neg}}

\newcommand {\Idm }{\mathrm{I}}
\newcommand {\Cent }{\mathrm{Cent}}
\newcommand {\disc }{\operatorname{disc}}

\newcommand {\Wind }{\operatorname{ind}}

\usepackage{ifdraft}
\ifdraft{\usepackage[outer]{showlabels}
	\usepackage{todonotes} \usepackage{datetime}}{}
\usepackage[pagebackref,colorlinks,citecolor=red,urlcolor=blue,linkcolor=green, bookmarks=false,hypertexnames=true]{hyperref}
\usepackage{orcidlink}
\usepackage{fancyhdr}
\usepackage{verbatim}

\begin{document}

\bibliographystyle{plain}

\setcounter{page}{1}

\thispagestyle{empty}


\author{Klaus Nielsen}\,\orcidlink{0009-0002-7676-2944}
\email{klaus@nielsen-kiel.de}
\ifdraft{\today \ \currenttime}{\date{November 15, 2024}}

\pagestyle{fancy}
\fancyhf{}
\fancyhead[OC]{Klaus Nielsen}
\fancyhead[EC]{Bireflectionality in real orthogonal groups}
%\fancyfoot[C]{\thepage}
\fancyhead[OR]{\thepage}
\fancyhead[EL]{\thepage}

\begin{abstract}
Let $\Omega(p,q)$ be the commutator subgroup of the orthogonal group $\OG(p,q)$ of signature $(p,q)$ over the reals. It is shown that an element of $\Omega(p,q)$ is bireflectional (product of 2 involutions) if and only if it is reversible (conjugate to its inverse). Moreover, the bireflectional elements 
of $\Omega(p, q)$ are classified.
\end{abstract}

\maketitle
%%%%%%%%%%%%%%%%%%%%%%%%%%%%%%%%%%%%%%%

%%%%%%%%%%%%%%%%%%%%%%%%%%%%%%%%%%%%%%%%%%%%%%%%%%%%%%%%%%%%%
\section{Introduction} \label{intro-sec}
We call an element of a group $G$ bireflectional if it is a product of 2 involutions of $G$.
We say that an element of $G$ is reversible if it is conjugate to its inverse. Clearly, a bireflectional element is reversible.
A group  is bireflectional if all its elements are bireflectional.
 
We are interested in classifying reversible and bireflectional elements in classical groups. 
It is well known that the  orthogonal group $\OG(V, Q)$ of a nondefective quadratic space $(V, Q)$ is bireflectional by theorems of Wonenburger \cite{Wonenburger-1966} and Ellers and Nolte \cite{EllersNolte-1982}. Wonenburger \cite{Wonenburger-1966} and \DJ okovi\'c \cite{Dokovic-1967} have shown that all reversible elements of the general linear group are bireflectional. 
In \cite{KN-1987}, F. Knüppel and the author classified bireflectional elements in special orthogonal groups.
Bireflectionality in unitary groups have been studied by Bünger \cite{Bunger-1997}, \DJ okovi\'c \cite{Dokovic-1971}, Gates et al. \cite{GSV-2014}, and Schaeffer Fry and Vinroot \cite{SFV-2016}. In a forthcoming paper, we deal with symplectic groups.

In \cite{KN-2010}, F. Knüppel and the author considered bireflectionality in
the commutator group $\Omega(p,q)$ of the orthogonal group of a symmetric bilinear space of signature $(p,q)$ over the reals.

\begin{theorem} \label{KN}
The following are are equivalent
\begin{enumerate}
	\item $\Omega(p,q)$ is bireflectional, 
	\item  $p, q, p+q \not \equiv2 \mod 4$, 
	\item All elements of $\Omega(p, q)$ are reversible.
\end{enumerate}
\end{theorem}

%%%%%%%%%%%%%%%%%%%%%%%%%%%%%%%%%%%%%%%
\section{Main results}

Let $V$ be a finite dimensional symmetric bilinear  space over the reals with signature $(p,q)$. Let $\Omega(p,q) = \Omega(V)$ denote the commutator subgroup of the orthogonal group $\OG(p,q) = \OG(V)$.
We show that all  reversible elements of $\Omega(p,q)$ are bireflectional and
give a complete classification of bireflectional elements of $\Omega(p, q)$. 

We need a definition:
For $\varphi \in \End(V)$ let
$\Bahn^{\infty}(\varphi^2)$ be the Fitting one space of $\varphi^2-1$ and
$n_j(\varphi)$ denote the number of elementary divisors $(x\pm1)^d$ of $\varphi$ of degree $d \equiv j \mod 8$. For a nondegenerate subspace $T$ of $V$ let
$q(T)$ denote the negative signature of $T$ and $p(T)$ denote the positive signature of $T$.

\begin{theorem} \label{main-th-1}
	Let $\varphi \in \Omega(p,q)$. Then $\varphi$ is bireflectional if and only if
	$\varphi$ is bireflectional in the special orthogonal group $\SOG(p,q)$ and one of the following holds
	\begin{enumerate}
		 \item $\varphi$ has an orthogonal summand of even dimension and  discriminant -1;
		 \item $2 n_3(\varphi) + 2 n_5(\varphi) + n_2(\varphi) + n_6(\varphi)
		 \equiv q(\Bahn^{\infty}(\varphi^2)) \mod 4$;
		 \item  $2 n_3(\varphi) + 2 n_5(\varphi) + n_2(\varphi) + n_6(\varphi)
		 \equiv p(\Bahn^{\infty}(\varphi^2))  \mod 4$.
	\end{enumerate}
\end{theorem}

\begin{theorem} \label{main-th-2}
	All reversible elements of $\Omega(p,q)$ are bireflectional. 
\end{theorem}

%%%%%%%%%%%%%%%%%%%%%%%%%%%%%%%%%%%%%%%
\section{Preliminaries}
 
 %%%%%%%%%%%%%%%%%%%%%%%%%%%%%%%%%%%%%%
 \subsection{Notations}

 For a linear mapping $\varphi$ of $V$ let $\Bahn^j(\varphi)$ denote the image and $\Fix^j(\varphi)$ the kernel of $(\varphi -1)^j$. Put $\Fix^{\infty}(\varphi) = \bigcup_{j\ge 1} \Fix^j(\varphi)$ and $\Bahn^{\infty}(\varphi) = \bigcap_{j\ge 1}
 \Bahn^j(\varphi)$. The space $\Bahn(\varphi) := \Bahn^1(\varphi)$ is the path or residual space of $\varphi$.
 And $\Fix(\varphi) := \Fix^1(\varphi)$ is the fix space of $\varphi$. Further, let $\Neg^j(\varphi) =  \Fix^j(-\varphi)$. The space $\Neg(\varphi) = \Neg^1(\varphi)$ is the negative space of  $\varphi$.
 
 For a monic polynomial $f$ of degree $d$ with $f(0) \ne 0$, let $f^*(x) = f(0)^{-1}x^d f(x^{-1})$ denote the reciprocal of $f$.
 By $\mu(\varphi)$ we denote the minimal polynomial of 
 $\varphi$. If $\varphi \in \OG(V)$, then $\mu(\varphi)$ is selfreciprocal; i.e. 
 $\mu(\varphi) = \mu(\varphi)^*$. 
 
 For a subspace $T$ of $V$, let $\disc T$ denote the discriminant of $T$.
 By $\Theta(\varphi)$, we denote the spinor norm or spinorial norm of an orthogonal transformation $\varphi$, as defined by Lipschitz or Eichler. See Artin \cite[Definition 5.5]{EArtin}.
Other authors like Dieudonn\'e or Wall use a slightly different definition;
see e.g. \cite[ch. 9,  3.4 Definition]{WScharlau}. But when restricted to the special orthogonal group, all spinor norms are equal.
The Witt index $\Wind V$ of $V$ is the dimension of a  maximal totally isotropic subspace of $V$. 

 We collect some well-known facts.
 
 \begin{remark} \label{remark-1}
 	We have $\det \varphi = \det \varphi|_{\Bahn(\varphi)}$.
 	If $\varphi \in \OG(V)$, then 
 	\begin{enumerate}
 		\item $\det \varphi = \det \varphi|_{\Neg(\varphi)}$.
 		\item $\varphi$ is involutory if and only if $\Bahn(\varphi)  = \Neg(\varphi)$.
 		\item $V f(\varphi)^{\perp} = \ker f^*(\varphi)$ for all polynomials $f$ with $f(0) \ne 0$.
 		\item $\Bahn^j(\varphi)^{\perp} = \Fix^j(\varphi)$.
 		\item $V = \Bahn^{\infty}(\varphi) \oplus \Fix^{\infty}(\varphi)$.
 		\item If $\sigma \in  \OG(V)$ is involutory, then $\Theta(\sigma) = \disc \Bahn(\sigma)$.
 		\item $\Omega(V) = \SOG(V) \cap \ker \Theta$.
 	\end{enumerate}	
 \end{remark}

 %%%%%%%%%%%%%%%%%%%%%%%%%%%%%%%%%%%%%%%
 \subsection{Orthogonally indecomposable transformations}
 
 Let $\OG(V)$ be the orthogonal group of a nondegenerate bilinear space $V$ over a field of characteristic not 2.
 
According to Huppert \cite[1.7 Satz]{Huppert-1980a}, an orthogonally indecomposable transformation of $\OG(V)$ is either\footnote{Our type enumeration follows Huppert \cite{Huppert-1980a}. In \cite{Huppert-1990} Huppert uses  a different  enumeration.}
\begin{enumerate}
	\item bicyclic with elementary divisors $e_1 = e_2 =(x \pm 1)^m$
	(type 1) or
	\item indecomposable as a linear transformation (type 2) or
	\item cyclic with minimal polynomial $(hh^*)^t$, where $h$ is irreducible and prime to its reciprocal $h^*$ ( type 3).
\end{enumerate}

In \cite{KN-1987a}, the authors proved the following

\begin{lemma}      
	                          \label{lemma1}
	 Let $\varphi \in \OG(V)$.  Assume that
	$\varphi$ is reversed by an involution $\sigma \in \OG(V)$.
	Then $V$ has a decomposition $V = U_1 \perp U_2 \perp \dots \perp U_m$, where the subspaces $U_j$ are $\varphi$- and $\sigma$-invariant and orthogonally indecomposable w.r.t. $\varphi$.
\end{lemma}

\begin{proof}
	See  \cite[Proposition p. 212]{KN-1987a} or more general, \cite[Proposition 6.1]{GKN-2008}.
\end{proof}

This reduces the problem of computing the determinant and spinor norm of a reversing involution to orthogonally indecomposable transformations.
It is convenient to introduce subtypes for the types 2:
\begin{definition}
	Let $\varphi \in \OG(V)$ be orthogonally indecomposable of type 2 with 
	$\mu(\varphi) = p^m$, where $p$ is irreducible.
	We say that $\varphi$ is of type 
	\begin{enumerate}
		\item $2\pm$ if $p = x \pm  1$, 
		\item $2*$ if $p \ne x \pm  1$,
		\item 2o if $p \ne x \pm 1$ and $m$ is  odd,
		\item 2e if $p \ne x \pm 1$ and  $m$ is even.
	\end{enumerate}	
\end{definition}

In \cite[2.3 Satz]{Huppert-1980a}, \cite[2.4 Satz]{Huppert-1980a}, and \cite[3.2 Satz]{Huppert-1980b}, Huppert proved the following:

\begin{proposition} \label{prop-1}
	Let $\ch K \ne 2$.
	Let $\varphi \in \OG(V)$ be orthogonally indecomposable. Then 
	\begin{enumerate}
		\item If $\varphi$ is of type 1, then $\dim V \equiv 0 \mod 4$.
		\item $\dim V$ is odd if and only if $\varphi$ is of type $2\pm$.
	\end{enumerate}
\end{proposition}

\begin{definition}
Let $\varphi \in \OG(V)$. We say that $\varphi$ is
of type T, if all orthogonally indecomposable  orthogonal  summands of $\varphi$ are of type T.
\end{definition}

\begin{lemma}      \label{lemma2}
	Let $\ch K \ne 2$.
	Let $\varphi \in \OG(V)$ be cyclic with minimal polynomial
	$(x-1)^n$.  Then $\varphi \in \Omega(V)$. Let $\varphi = \sigma \tau$ be a product of two orthogonal involutions. Then 
	\begin{enumerate}
		\item $\dim \Bahn(\sigma) = \dim \Bahn(\tau) \in \{m,m+1\}$, where $m = \lfloor \frac{n}{2} \rfloor$.
		\item Either $\Bahn(\sigma)$ and $\Bahn(\tau)$ or
		$\Fix(\sigma)$ and $\Fix(\tau)$ are hyperbolic spaces.
		\item $\Theta(\sigma) = \Theta(\tau) \in \{(-1)^t, (-1)^t \disc V \}$,
		where $t = \lfloor \frac{n+1}{4} \rfloor = \lfloor \frac{m+1}{2} \rfloor$.
	\end{enumerate}
\end{lemma}

\begin{proof}
	1 is clear and is shown in \cite[Lemma 2.6]{KhN-1987}.
	The subspace $T = \Bahn^m(\varphi)$ is totally isotropic of dimension $m$ and
	$\langle \sigma, \tau \rangle$-invariant.
	Then  $ T = [T \cap  \Bahn(\sigma)] \oplus [T \cap  \Fix(\sigma)]$
	and  $ T = [T \cap  \Bahn(\tau)] \oplus [T \cap  \Fix(\tau)]$.
	We may assume that $\dim \Bahn(\sigma) = m = \dim \Bahn(\tau)$.
	If $m$ is even, then
	$\dim [T \cap  \Bahn(\tau)] = m/2 = \dim [T \cap  \Bahn(\sigma)]$
	and $\Bahn(\tau), \Bahn(\sigma)$ are hyperbolic.
	If $m$ is odd, then $\dim [T \cap \Fix(\tau)] = \frac{1}{2} (m+1) =
	\dim [T \cap \Fix(\sigma)]$
	and $\Fix(\sigma)$ and $\Fix(\tau)$ are hyperbolic.
\end{proof}

\begin{lemma}                                                                   \label{lemma3}
	Let $V$ be hyperbolic, and let $T$ be a totally isotropic subspace of $V$ with $\dim T = \frac{\dim V}{2}$.
Let $\varphi \in \OG(V)$, 
and let $\sigma \in \OG(V)$ be an involution reversing $\varphi$. If $T$ is $\sigma$-invariant, 
then $\Bahn(\sigma)$ and  $\Fix(\sigma)$ are hyperbolic spaces
of dimension $\frac{1}{2} \dim V$.
\end{lemma}

\begin{corollary}                                                  \label{lemma4}
	Let $\varphi \in \OG(V)$ be orthogonally indecomposable of 1.
	Let $\sigma \in \OG(V)$ be an involution reversing $\varphi$.
	Then
	\begin{enumerate}
		\item $\Bahn(\sigma)$ and  $\Fix(\sigma)$ are hyperbolic.
		\item $\dim \Bahn(\sigma) = \dim \Fix(\sigma) = \frac{1}{2} \dim V$.
		\item $\varphi \in \Omega(V)$.
		\item $\det \sigma = 1$, and $\Theta(\sigma) = (-1)^{\frac{1}{4} \dim V}$.
		\end{enumerate}
\end{corollary}

\begin{lemma}                                                  \label{lemma5}
Let $\varphi \in \OG(p,q)$ be orthogonally indecomposable of type $2^-$.
Then $p-q = \pm 1$, $\varphi \in \Omega(p,q)$.
Let $\sigma \in \OG(p,q)$ be an involution reversing $\varphi$. Then
\begin{enumerate}
	\item If $n = 8t +1$, then $q \in \{4t, 4t+1\}$ and $\Theta(\sigma) = \det \sigma^q$.
	\item If $n = 8t +7$, then $q \in \{4t+3, 4t+4\}$ and $\Theta(\sigma) = \det \sigma^q$.
	\item If  $n = 8t +3$, then $q \in \{4t+1, 4t+2\}$ and $\Theta(\sigma) = -\det \sigma^q$.
	\item If  $n = 8t +5$, then $q \in \{4t+2, 4t+3\}$ and $\Theta(\sigma) = -\det \sigma^q$.
\end{enumerate}
If $q=2r$ is even, then $\Theta(\sigma) = (-1)^r$.
\end{lemma}

\begin{proof}
See \cite[Lemma 5.1]{KT-1998}.
\end{proof}

\begin{lemma}                                                  \label{lemma6}
	Let $\varphi \in \OG(p,q)$ be orthogonally indecomposable
	of type  $2^*$.  Then $q=2r$ is even and $\varphi \in \Omega(p,q)$.
	Let $\sigma \in \OG(p,q)$ be an involution reversing $\varphi$. Then
	$\Theta(\sigma) = (-1)^r$ and
	\begin{enumerate}
		\item If $\varphi$ is of type 2e, then  $q = p$ and $\det \sigma = 1$.
		\item If $\varphi$ is of type 2o, then  $|q - p| = 2$ and $\det \sigma = -1$.
	\end{enumerate}
\end{lemma}

\begin{proof}
	See also \cite[Lemma 5.3]{KT-1998} and \cite[Lemma 5.4]{KT-1998}.
	If $\varphi$ is of type 2e apply \ref{lemma3}. Let $\varphi$ be of type 2o, and let $\mu(\varphi) = g^{2t+1}$, where $g$ is irreducible (of degree 2). Then $\ker g(\varphi)^t = V g(\varphi)^{t+1}$ is totally isotropic and $\sigma$-invariant. Hence $\Wind \Bahn(\sigma) = \Wind \Fix(\sigma) = t$.
	By \cite[4.1 Satz]{Huppert-1980b}, $\disc V = 1$ so that $q = 2t$ or $q = 2t+2$. Further, $q(\Bahn(\sigma)) = \frac{q}{2}$.
\end{proof}

\begin{lemma}                                                  \label{lemma7}
	Let $\varphi \in \OG(p,q)$ be orthogonally indecomposable of type $3$
	with minimum polynomial $(x-\lambda)^t (x-\lambda^{-1})^t$.
	Let $\sigma \in \OG(p,q)$ be an involution reversing $\varphi$.
	Then $p=q=t$ and $\Theta(\varphi) = \lambda^q$
	and
	\begin{enumerate}
		\item If $q=2r$ is even, then $\Theta(\sigma) = (-1)^r$.
		\item If $q$ is odd, then $\Theta(-\sigma) = -\Theta(\sigma)$.
	\end{enumerate}
\end{lemma}

\begin{lemma}                                                  \label{lemma8}
	Let $\varphi \in \OG(p,q)$ be orthogonally indecomposable of type 3
	with minimum polynomial $(h h^*)^t$, where $h$ is irreducible of degree 2.
	Then $q=p=2t$, and $\Theta(\varphi) = h(0)^t > 0$.
	Furthermore, there exists an orthogonal involution $\sigma$ reversing
	$\varphi$ such that $\Theta(\sigma) = (-1)^t$.
\end{lemma}

\begin{proof}
	 In a suitable basis, we have
	\[
	\varphi = \left (\begin{array} {cc} A & 0\\ 0 & A^+ \end{array} \right ),
	f = \left (\begin{array} {cc} 0 & \Idm_p\\ \Idm_p & 0\end{array} \right ),
	\]
	where $A^+$ is the transpose inverse of $A$.
	It is well known that $A$ is transposed by a symmetric matrix $S$ with maximal Witt index ($= t$); cf. e.g. \cite[Theorem 66]{Kaplansky-1974}.
	Let
	\[
	\sigma = \left (\begin{array} {cc} 0 & S\\ S^{-1} & 0\end{array} \right ).
	\]
	
	Then $\sigma \in \OG(p, q)$, $\varphi^{\sigma} = \varphi^{-1}$, and $\Theta(\sigma) = \det S = (-1)^t$.	
\end{proof}

\begin{lemma}                        \label{lemma9}
	Let $\varphi \in \OG(p,q)$.
	Assume that  every orthogonally indecomposable summand of $\varphi$
	has dimension $\equiv 0 \mod 4$ or is of type 2*.
	Then $\Cent(\varphi) \subseteq  \Omega(p,q)$.
\end{lemma}

\begin{proof}
	Clearly, $n$ is even. Every orthogonally indecomposable summand of $\varphi$ is either hyperbolic or of type 2o.
	Using \ref{lemma6}, we see that $\disc V = 1$.
	Let $\xi \in \Cent(\varphi)$.
	We may assume that $\xi$ is orthogonally primary.

	If $\xi$ is unipotent, then clearly $\Theta(\xi)=1 = \det \xi$.
	If $-\xi$ is unipotent, then $\Theta(\xi)=\Theta(-\xi) \disc V = 1$ and $\det \xi = (-1)^n \det (-\xi) = 1$.
	If $\xi$ is of type $3$ and split, then $\Theta(\xi) > 0$ by \ref{lemma7}.
	If $\xi$ is of type $3$ and nonsplit, then $\Theta(\xi) > 0$ by \ref{lemma8}.
	If $\xi$ is of type $2^*$, then $\xi \in \Omega(p,q)$ by \ref{lemma6}.
\end{proof}

\begin{corollary}                                              \label{lemma10}
	Let $\varphi \in \OG(p,q)$ be orthogonally indecomposable of type 3
	with minimum polynomial $(h h^*)^t$, where $h$ is irreducible of degree 2.
	Then $\Theta(\sigma) = (-1)^t$ for all orthogonal involution $\sigma \in \OG(p,q)$ reversing $\varphi$.
\end{corollary}

%**********************************************************************
%**********************************************************************
\section{Proof of theorem \ref{main-th-1}}

\begin{lemma}      \label{lemma11}
	Let $\ch K \ne 2$.
	Let $\varphi \in \OG(V)$ be orthogonally indecomposable of type 2-.
	Assume that $\dim V \ge 3$. Then $\disc V = -\disc [\Bahn(\varphi)/ \Fix(\varphi)]$.
\end{lemma}

\begin{proof}
	Let $T$ be a complement of $\Fix(\varphi)$ in $\Bahn(\varphi)$. Then
	$T^{\perp}$ is a hyperbolic plane, and $T$ is isomorphic to $\Bahn(\varphi)/ \Fix(\varphi)$. Then $\disc V = (\disc T) (\disc T^{\perp}) = - \disc T = - \disc [\Bahn(\varphi)/ \Fix(\varphi)]$.
\end{proof}	

Using Sylverster's law of inertia, we obtain

\begin{corollary} \label{cor-x}
	Let $\varphi \in \OG(p,q)$ be of type $2\pm$.  Assume  that  $\varphi$ has orthogonal  summands $\varphi_1$ and $\varphi_2$ of odd dimension with $\disc \varphi_1 = 1$ and $\disc \varphi_2 = -1$. Then $\varphi$ has an orthogonal  summand of even dimension and discriminant -1.
\end{corollary}

\begin{lemma} \label{lemma12}
	Let $\varphi \in \OG(p,q)$ be of type $2\pm$.  Assume further that all orthogonally indecomposable orthogonal  summands of $\varphi$ have discriminant $\delta$. Let $\sigma \in \OG(p,q)$ be an involution reversing $\varphi$. Then 
	\begin{enumerate}
		\item If $\delta = 1$, then $\Theta(\sigma) = (-1)^{n_3(\varphi) + n_5(\varphi)}$.
		\item If $\delta = -1$, then $\Theta(\sigma) = (-1)^{n_3(\varphi) + n_5(\varphi)} \det \sigma$.
	\end{enumerate}
\end{lemma}

\begin{proof}
	By \ref{lemma1}, we may assume that $\varphi$ is orthogonally indecomposable.
	Apply \ref{lemma5}.
\end{proof}	

\begin{proof}[Proof of theorem \ref{main-th-1}]
	If  $\varphi \in \Omega(p, q)$ has an orthogonal  summand of type $2\pm$ with even dimension and discriminant -1, then clearly $\varphi$ is reversed by an involution with prescribed spinor norm and prescribed determinant. 
	
	Next assume that $\varphi$ has an orthogonal summand of type 3 with discriminant -1.
	Then $\varphi$ is reversed by an involution $\rho \in \OG(p,q)$ with $\Theta(\rho) = 1$. 
	If $\varphi$ has an orthogonally summand of type $2\pm$ we can additionally achieve that $\det \rho = 1$. 
	If $\varphi$ has no orthogonally summand of type $2\pm$, then $\det \rho = 1$ as $\varphi$ is bireflectional in $\SOG(p,q)$ and all reversing involutions have the same determinant $(-1)^{\frac{\dim V}{2}}$.
	
	So we may assume that $\varphi$ has no orthogonal summand of type 3 with discriminant -1.
	And by \ref{cor-x}, we may assume that all orthogonally indecomposable orthogonal summands of type $2\pm$ of $\varphi$ have the same discriminant
	 $\delta$.
	 
	Let $\sigma \in \OG(p,q)$ be an involution reversing $\varphi$.
	Put $B = \Bahn^{\infty}(\varphi^2), F = \Fix^{\infty}(\varphi^2)$.
	Assume first that $\delta = 1$. 
	Then $\Theta(\sigma) = 1$ iff 
	$n_3(\varphi) +  n_5(\varphi) + \frac{n_2(\varphi)}{2} + \frac{n_6(\varphi)}{2} \equiv \frac{q(B)}{2} \mod 2$.
	If $\varphi$ has no orthogonal summand of type $2\pm$, then $\dim V \equiv 0 \mod 4$ so that $\det \sigma = 1$. Otherwise,  we can adjust the determinant of $\sigma$ without changing its spinor norm.
	
	Finally let $\delta = -1$.
	Then $\sigma \in \Omega(p,q)$ if and only if
	\begin{enumerate}
		\item $\det \sigma_F = \det \sigma_B = (-1)^{\frac{\dim B}{2}}$ and
		\item $(-1)^{\frac{q(B)}{2}} = \Theta(\sigma_B) = \Theta(\sigma_F) =   (-1)^{\frac{n_2(\varphi)}{2} + \frac{n_6(\varphi)}{2}}(-1)^{n_3(\varphi) +  n_5(\varphi)} \det \sigma_F$.
	\end{enumerate}
\end{proof}

%**********************************************************************
%**********************************************************************
\section{Proof of theorem \ref{main-th-2}}

\begin{lemma}                        \label{lemma13}
	Let $\varphi \in \OG(p,q)$.
	Assume  that  
	all  orthogonally indecomposable summands of type $2\pm$ of $\varphi$  have  discriminant -1 and all orthogonally indecomposable summands of type $3$ of $\varphi$ have discriminant 1.
	Then $\Theta(\xi) = \det \xi $ for all $\xi \in \Cent(\varphi)$. 
\end{lemma}

\begin{proof}
	We may assume that $\xi$ is orthogonally primary.
	If $\xi$ is unipotent, then  $\xi \in \Omega(p,q)$ by \ref{lemma2} and \ref{lemma4}.
	If $-\xi$ is unipotent, then clearly $\Theta(\xi) = \disc V = (-1)^{\dim V} = \det \xi$. If $\xi$ is of type 2* or 3(nonsplit), then
	$\xi \in \Omega(p,q)$ by \ref{lemma6} and \ref{lemma8}.
	So let $\xi$ be of split type 3. Then $\disc V = 1$ as $\dim V$ is even. Hence $\dim V \equiv 0 \mod 4$ as $V$ is hyperbolic.  By \ref{lemma7}, $\Theta(\xi) = 1$.
\end{proof}

Similarly, we have 
\begin{lemma}                        \label{lemma14}
	Let $\varphi \in \OG(p,q)$. 
	Assume  that  all  orthogonal summands of $\varphi$  have discriminant +1.
	Then $\Cent(\varphi) \subseteq \ker \Theta$. 
\end{lemma}

\begin{proof}[Proof of theorem \ref{main-th-2}]
	Let $\varphi, \alpha \in \Omega(p,q)$ with $\varphi^{\alpha} = \varphi^{-1}$. 
	By \cite[Theorem 4.1]{KN-2010}, $\varphi$ is reversed by an involution $\sigma \in \SOG(p,q)$. Then $\alpha \sigma$ commutes with $\varphi$. 
	By \ref{lemma14}, we are already done if $\varphi$ has no orthogonal summand of discriminant -1. Further by \ref{main-th-1}, we may assume that  $\varphi$ has no orthogonal summand of even dimension  with discriminant -1.
	So all orthogonal summands of $\varphi$ with odd dimension  have discriminant -1. By \ref{lemma13}, $\Theta(\alpha \sigma) = 1$.
\end{proof}

%%%%%%%%%%%%%%%%%%%%%%%%%%%%%%%%%%%%%%%%%%%%%%%%%%%%%%%%%%%%%

\end{document}

\typeout{get arXiv to do 4 passes: Label(s) may have changed. Rerun}